\documentclass[12pt]{amsart}
\usepackage{latexsym,amssymb,amscd}

\def\B'c{{\mathcal{B'}}}
\def\U'c{{\mathcal{U'}}}

\def\opn#1#2{\def#1{\operatorname{#2}}} 
\opn\chara{char}
\opn\length{\ell}
\opn\projdim{proj\,dim}
\opn\injdim{inj\,dim}
\opn\ini{in}
\opn\rank{rank}
\opn\Tiefe{Tiefe}
\opn\grade{grade}
\opn\height{height}
\opn\embdim{emb\,dim}
\opn\codim{codim}

\opn\Tr{Tr}
\opn\bigrank{big\,rank}
\opn\superheight{superheight}\opn\lcm{lcm}
\opn\trdeg{tr\,deg}%
\opn\reg{reg}
\opn\lreg{lreg}
\opn\deg{deg}
\opn\lcm{lcm}
\opn\div{div}
\opn\Div{Div}
\opn\cl{cl}
\opn\Cl{Cl}
\opn\Spec{Spec}
\opn\Supp{Supp}
\opn\supp{supp}
\opn\Sing{Sing}
\opn\Ass{Ass}
\opn\Ann{Ann}
\opn\Rad{Rad}
\opn\Soc{Soc}
\opn\Ker{Ker}
\opn\Coker{Coker}
\opn\Im{Im}
\opn\Hom{Hom}
\opn\Tor{Tor}
\opn\Ext{Ext}
\opn\End{End}
\opn\Aut{Aut}
\opn\id{id}

\opn\nat{nat}
\opn\GL{GL}
\opn\SL{SL}
\opn\mod{mod}
\opn\ord{ord}
\opn\depth{depth}
\opn\aff{aff}
\opn\con{conv}
\opn\relint{relint}
\opn\st{st}
\opn\lk{lk}
\opn\cn{cn}
\opn\core{core}
\opn\vol{vol}
\opn\gr{gr}

\def\pot#1#2{#1[\kern-0.28ex[#2]\kern-0.28ex]}

\opn\dirlim{\underrightarrow{\lim}}
\opn\invlim{\underleftarrow{\lim}}
\let\union=\cup

\let\dirsum=\oplus

\def\pnt{{\raise0.5mm\hbox{\large\bf.}}}

\def\Implies{\ifmmode\Longrightarrow \else
     \unskip${}\Longrightarrow{}$\ignorespaces\fi}
\def\implies{\ifmmode\Rightarrow \else
     \unskip${}\Rightarrow{}$\ignorespaces\fi}
\def\iff{\ifmmode\Longleftrightarrow \else
     \unskip${}\Longleftrightarrow{}$\ignorespaces\fi}

\let\:=\colon
\newtheorem{Theorem}{Theorem}[section]
\newtheorem{Lemma}[Theorem]{Lemma}
\newtheorem{Corollary}[Theorem]{Corollary}
\newtheorem{Proposition}[Theorem]{Proposition}

\newtheorem{Example}[Theorem]{Example}

\newtheorem{Definition}[Theorem]{Definition}

\let\epsilon=\varepsilon
\let\phi=\varphi
\let\kappa=\varkappa
\title{Constructible ideals}
\author{}

\begin{document}
\flushbottom
\maketitle

\begin{center}

Anda Olteanu\footnote{Work supported by the CEEX Programme of the Romanian Ministry of Education and Research, contract CEX 05-D11-11/2005 and by the CNCSIS contract TD 101/2007.}
\end{center}
	\[
\]

\section*{Introduction}

Properties of simplicial complexes such as Cohen-Macaulayness or shellability can be obtained by studying the Stanley-Reisner ideal of the Alexander dual. More precisely, by Eagon-Reiner theorem \cite{EaRe} a simplicial complex, $\Delta$, is Cohen-Macaulay if and only if the Stanley-Reisner ideal of the Alexander dual, $I_{\Delta^{\vee}}$, has a linear resolution. Also, a pure simplicial complex is shellable if and only if the Stanley-Reisner ideal of the Alexander dual has linear quotients \cite{HeHiZh}.

For pure simplicial complexes, the following implications are known: 
\[\mbox{shellable}\Rightarrow\mbox{constructible}\Rightarrow\mbox{Cohen-Macaulay}.
\]

A natural question arises: Determine the Stanley-Reisner ideal of the Alexander dual of a constructible simplicial complex. Our paper aims to answer this question.

The paper is structured as follows. In Section 1 we recall the notion of simplicial complex and some results related to this concept.

In Section $2$ we introduce the concept of constructible ideal and we relate this concept with the notion of constructible simplicial complex.

Next, in Section $3$ we prove that every constructible ideal has a linear resolution and we find a formula for the Betti numbers of a constructible ideal.

In Section $4$ we show that the polarization of a constructible ideal is a square-free constructible ideal.

In Section $5$ we find some properties of monomial ideals with linear quotients and we prove that every monomial ideal with linear quotients is a constructible ideal.

Finally, in Section $6$, we discuss some examples.\\

The author would like to thank Professor J\"urgen Herzog for valuable suggestions and comments during the preparation of this paper.

\section{Basic facts} First, we recall the notion of simplicial complex and some concepts related to it. We denote by $S=k[x_1,\ldots, x_n]$ the polynomial ring in $n$ variables over a field $k$.

\begin{Definition}\rm
 Let $n>0$ be an integer and $[n]=\{1,\ldots, n\}$. A \it simplicial complex, \rm$\ $ $\Delta$, on $[n]$ is a colection of subsets of $[n]$ that satisfy the following conditions:
\begin{itemize}
	\item [(a)] $\{i\}\in \Delta$, for all $i\in [n]$;
	\item[(b)] If $F\in \Delta$ and $G$ is a subset of $F$, then $G\in\Delta$.
\end{itemize}
\end{Definition}

$[n]$ is called \it the vertex set \rm $\ $of $\Delta $ and the elements of the simplicial complex are called \it faces \rm. The \it dimension \rm of a face, $F$, is denoted $\dim(F)$ and $\dim(F)=|F|-1$. Denote by $d=\max\{|F|\ :\ F\in \Delta\}$. Then the dimension of the simplicial complex $\Delta$ is $d-1$. A \it facet\rm$\it$ of $\Delta$ is a maximal face (with respect to the inclusion). If all the facets have the same dimension, we say that the simplicial complex is \it pure\rm. 

Denote by $\mathcal{F}(\Delta)$ the set of all the facets of $\Delta$.

Let $\Delta$ be a simplicial complex on the vertex set $[n]$. Define the ideal
	\[I_{\Delta}=(x_F\ :\ F\notin\Delta)\subset S,
\]
where $x_F=x_{i_1}\ldots x_{i_s}$ for $F=\{i_1,\ldots,i_s\}$. $I_{\Delta}$ is called the \it Stanley-Reisner ideal \rm $\ $of $\Delta$. Also, it can be considered \it the Stanley-Reisner ring\rm
	\[k[\Delta]=\frac{k[x_1,\ldots,x_n]}{I_{\Delta}}.
\]

\it The facet ideal\rm$\ $ of $\Delta$, $I(\Delta)$, is the square-free monomial ideal generated by all the monomials $x_F=\prod\limits_{i\in F}x_i$, where $F$ is a facet of $\Delta$. So, if $\Delta=\langle F_1,\ldots,F_r\rangle$,
	\[I(\Delta)=(x_{F_1},\ldots,x_{F_r}).
\]

\begin{Definition}\rm
 Let $\Delta$ be a simplicial complex on the vertex set $[n]$. The simplicial complex
	\[\Delta^{\vee}=\{[n]\setminus F\ :\ F\notin\Delta\}
\]
is called \it the Alexander dual\rm  $\ $of $\Delta$.
\end{Definition}

For the simplicial complex $\Delta$ with the vertex set $[n]$, let $\Delta^c$ be the simplicial complex whose facets are $[n]\setminus F$, where $F$ is a facet in $\Delta$.

If $\Delta$ is a simplicial complex, there exists the following relation between $\Delta^{\vee}$ and $\Delta^c$:

\begin{Proposition}\cite{HeHiZh}
\it Let $\Delta $ be a simplicial complex. Then
	\[I_{\Delta^{\vee}}=I(\Delta^c).
\]
\end{Proposition}

We recall the definition of a shellable simplicial complex.

\begin{Definition}\rm
 A pure simplicial complex, $\Delta$, is called \it shellable\rm$\ $ if its facets can be ordered $F_1,\ldots, F_m$ such that, for all $2\leq i\leq m$, the simplicial complex
	\[\langle F_1,\ldots,F_{i-1}\rangle\cap\langle F_i\rangle
\]
is generated by maximal proper subsets of $F_i$.
\end{Definition}

To describe the connection between a shellable simplicial complex and its Alexander dual we need to recall the definition of a monomial ideal with linear quotients.

\begin{Definition}\rm\cite{HeTa}
 A monomial ideal $I$ of $S$ with the minimal system of generators $G(I)=\{u_1,\ldots,u_r\}$ is called \it ideal with linear quotients \rm with respect to the sequence of monomials $u_1,\ldots, u_r$ if, for all $\ 2\leq i\leq r$ and for all $j<i$, there exist $l$ and $k$, $l\in\{1,\ldots,n\}$ and $k<i$, such that $u_k/[u_k,u_i]=x_l$ and $x_l$ divides $u_j/ [u_j,u_i]$ where we denote $[u_j,u_i]:=\gcd(u_j,u_i)$. 
\end{Definition}

The connection between a shellable simplicial complex and its Alexander dual is given by the following result:

\begin{Theorem}\cite{HeHiZh}
\it Let $k$ be a field and $\Delta$ be a pure simplicial complex. Then $\Delta$ is shellable if and only if $I_{\Delta^{\vee}}$ has linear quotients.
\end{Theorem}

We say that a simplicial complex $\Delta$ is Cohen-Macaulay if the Stanley-Reisner ring, $k[\Delta]$, is Cohen-Macaulay.

It is known that a shellable simplicial complex is Cohen-Macaulay over every field, \cite{Hi}. 

The connection between a Cohen-Macaulay simplicial complex and its Alexander dual is given by the following result:
\begin{Theorem}
\rm(Eagon-Reiner)\cite{EaRe}\it$\ $ Let $k$ be a field and $\Delta$ be a simplicial complex. Then $k[\Delta]$ is Cohen-Macaulay if and only if $I_{\Delta^{\vee}}$ has a linear resolution.
\end{Theorem}

\begin{Definition}\rm\cite{St} A pure simplicial complex $\Delta$ is \it constructible\rm $\ $ if it can be obtained by the following recursive procedure:
\begin{itemize}
	\item[(i)] Any simplex is constructible.
	\item[(ii)] If $\Delta_1$ and $\Delta_2$ are constructible simplicial complexes of the same dimension, $d$, and if $\Delta_1\cap\Delta_2$ is a constructible simplicial complex of dimension $d-1$, then $\Delta_1\cup\Delta_2$ is a constructible simplicial complex. 
\end{itemize}
\end{Definition}

It is known that all shellable simplicial complexes are constructible, \cite{St}. Indeed, if $\Delta$ is a shellable simplicial complex with the shelling order of the facets $F_1,\ldots, F_m$, then we can consider $\Delta_1=\langle F_1,\ldots, F_{m-1}\rangle$ and $\Delta_2=\langle F_m\rangle$.

\begin{Theorem}\cite{BrHe}
\it A constructible simplicial complex is Cohen-Macaulay over every field.
\end{Theorem}

So, for pure simplicial complexes, we have the following implications:
	\[\mbox{shellable}\Rightarrow\mbox{constructible}\Rightarrow\mbox{Cohen-Macaulay}.
\]
If we consider the connection with the Stanley-Reisner ideal of the Alexander dual for a pure simplicial complex, $\Delta$, we have the following diagram:
	\[\begin{array}{ccccc}
	\Delta\ \mbox{is shellable}&\Longrightarrow&\Delta\ \mbox{is constructible}&\Longrightarrow&\Delta\ \mbox{is Cohen-Macaulay}\\
	\Updownarrow&&&&\Updownarrow\\
	I_{\Delta^{\vee}}\ \mbox{is an ideal with}&&&&I_{\Delta^{\vee}}\ \mbox{is an ideal with a}\\
	\mbox{linear quotients}&&&&\mbox{linear resolution} 
\end{array}
\]
The following question arises: how to characterize the Stanley-Reisner ideal of the Alexander dual associated to a constructible simplicial complex. For this, we introduce the concept of constructible ideal.
\section{Constructible ideals}

Let $S=k[x_1,\ldots, x_n]$ be a polynomial ring in $n$ variables over a field $k$. 
\begin{Definition}\rm
A monomial ideal $I$ of $S$ generated in degree $q$ is a \it constructible ideal \rm if it can be obtained by the following recursive procedure:
\begin{itemize}
	\item [(i)] If $u$ is a monomial in $S$ and $I=(u)$, then $I$ is a constructible ideal;
	\item[(ii)] If $I_1$, $I_2$ are constructible ideals generated in degree $q$ and $I_1\cap I_2$ is a constructible ideal  generated in degree $q+1$, then $I_1+I_2$ is a  constructible ideal.
\end{itemize}
\end{Definition}

We note that the recursion procedure will stop because if the ideal $I$ has the minimal system of generators $G(I)=\{u_1,\ldots,u_r\}$ and, if we consider $I=I_1+I_2$, in $I_1\cap I_2$ the generators can contain each variable to a power which is less or equal to the maximal power to which that variable appears in all the generators of $I$. Let $a_i$ be the maximum of the exponents of the variable $x_i$ in the generators of $I$ and let $\mathbf{\underline{a}}=(a_1,\ldots, a_n)$. The recursion procedure will stop after at most $|\mathbf{\underline{a}}|:=a_1+\ldots +a_n$ steps.

The above remarks show that we could consider also the following definition of the constructible ideals.

Let $\mathbf{\underline{a}}=(a_1,\ldots, a_n)\in\mathbb{Z}_{>0}^n$. We denote $$\mathcal{M}_{\mathbf{\underline{a}}}=\{x_1^{b_1}\ldots x_n^{b_n}\ :\ 0\leq b_i\leq a_i\ \mbox{for all} \ 1\leq i\leq n\}$$ and $$\mathcal{J}_{\mathbf{\underline{a}}}=\{I\ :\ I\ \mbox{monomial ideal of}\ S\ \mbox{with}\ G(I)\subseteq\mathcal{M}_{\mathbf{\underline{a}}}\}.$$ By $G(I)$ we mean the minimal system of generators of the monomial ideal $I$. We set $|\mathbf{\underline{a}}|=a_1+\ldots+a_n$.

We note that, if $I,\ J\in\mathcal{J}_{\mathbf{\underline{a}}}$, then $I\cap J\in\mathcal{J}_{\mathbf{\underline{a}}}$.

\begin{Definition}\rm Let $I\in\mathcal{J}_{\mathbf{\underline{a}}}$ be a monomial ideal generated in degree $q$. $I$ is an \it{$\mathbf{\underline{a}}-$constructible ideal}\rm$\ $ if it can be obtained by the following recursive procedure:
\begin{itemize}
	\item [(i)] If $u\in\mathcal{M}_{\mathbf{\underline{a}}}$ and $I=(u)$, then $I$ is an $\mathbf{\underline{a}}-$constructible ideal;
	\item[(ii)] If $I_1,\ I_2\in\mathcal{J}_{\mathbf{\underline{a}}}$ are $\mathbf{\underline{a}}-$constructible ideals generated in degree $q<|\mathbf{\underline{a}}|$ and $I_1\cap I_2\in\mathcal{J}_{\mathbf{\underline{a}}}$ is an $\mathbf{\underline{a}}-$constructible ideal  generated in degree $q+1$, then $I_1+I_2$ is an $\mathbf{\underline{a}}-$constructible ideal.
\end{itemize}
\end{Definition}

Note that an $\mathbf{\underline{1}}-$constructible ideal is a square-free monomial ideal, where $\mathbf{\underline{1}}=(1,\ldots,1)\in\mathbb{Z}_{>0}^n$.

It is also important to notice that the only  principal ideal in $\mathcal{J}_{\mathbf{\underline{a}}}$ generated in
degree $|\mathbf{\underline{a}}|$ is $  I= (x_1^{a_1}\cdot\ldots \cdot x_n^{a_n}) $. This observation justifies that the recursion procedure of the above definition eventually terminates.

It is obvious that a monomial ideal $I$ is a constructible ideal(in the sense of Definition $2.1$) if and only if $I$ is an $\mathbf{\underline{a}}-$constructible ideal, for some $\mathbf{\underline{a}}\in\mathbb{Z}_{>0}^n$. Although Definition $2.2$ looks more technical, it will turn out that it is very useful in the proofs.
\begin{Theorem}
\it Let $\Delta$ be a pure simplicial complex on the vertex set $[n]=\{1,\ldots,n\}$. The following are equivalent:
\begin{itemize}
	\item [(a)] $\Delta$ is constructible.
	\item[(b)] $I_{\Delta^{\vee}}$ is a square-free constructible ideal .
\end{itemize}
\end{Theorem}
\begin{proof} Actually, we show that $\Delta$ is a constructible simplicial complex if and only if  $I_{\Delta^{\vee}}$ is an $\mathbf{\underline{1}}-$constructible ideal.  

$"(a)\Rightarrow (b)"$ We use induction on the dimension of $\Delta$. 

If $\dim(\Delta)=1$ the statement is obvious. In this case, $\Delta$ is a Cohen-Macaulay simplicial complex; thus is connected.

Let $\Delta$ be a $d-$dimensional constructible simplicial complex. We prove by induction on the number of facets of $\Delta$ that $I_{\Delta^{\vee}}$ is an $\mathbf{\underline{1}}-$constructible ideal.

If $\Delta$ is a simplex, $\Delta=\langle F\rangle$, we have that $I_{\Delta^{\vee}}=(x_{F^c})$ and it is an $\mathbf{\underline{1}}-$constructible ideal, by definition.

Let $\Delta$ be a $d-$dimensional constructible simplicial complex with the facet set $\mathcal{F}(\Delta)=\{F_1,\ldots, F_r\}$, $r\geq2$. Since $\Delta$ is constructible, there exist two $d-$dimensional constructible subcomplexes $\Delta_1$ and $\Delta_2$ such that $\Delta=\Delta_1\union\Delta_2$ and $\Delta_1\cap\Delta_2$ is a $(d-1)-$dimensional constructible simplicial complex. Since $\Delta_1$ and $\Delta_2$ are $d-$dimensional constructible simplicial complexes with at most $r-1$ facets, $I_{\Delta_1^{\vee}}$ and $ I_{\Delta_2^{\vee}}$ are $\mathbf{\underline{1}}-$constructible ideals generated in degree $n-d-1$. 

Since $I_{\Delta_1^{\vee}}\cap I_{\Delta_2^{\vee}}=I_{(\Delta_1\cap \Delta_2)^{\vee}}$ and $\Delta_1\cap \Delta_2$ is a $(d-1)-$dimensional constructible simplicial complex, by the induction hypothesis we have that $I_{\Delta_1^{\vee}}\cap I_{\Delta_2^{\vee}}$ is an $\mathbf{\underline{1}}-$construc- tible ideal generated in degree $n-d$ and hence $I_{\Delta^{\vee}}=I_{\Delta_1^{\vee}}+I_{\Delta_2^{\vee}}$ is an $\mathbf{\underline{1}}-$constructible ideal.

"$(b)\Rightarrow (a)$" We use descending induction on the degree of the monomials from the minimal system of generators of the monomial ideal $I$.

If $I=(x_1\ldots x_n)$, then $I=I_{\langle\{\emptyset\}\rangle}$ and $\langle\{\emptyset\}\rangle^{\vee}$ is the $n-$simplex, hence is a constructible simplicial complex.

Let $I$ be an $\mathbf{\underline{1}}-$constructible ideal generated in degree $q<n$. We use induction on the number of monomials from the minimal system of generators of the ideal $I$.

If $u\in\mathcal{M}_{\mathbf{\underline{1}}}$ and $I=(u)$, then let $\Delta$ be such that $I_{\Delta}=(u)$ and let $F=\supp(u)$. Then $\Delta^{\vee}$ is the simplex generated by $F^c$ and it is a constructible simplicial complex. 

Let $I$ be an $\mathbf{\underline{1}}-$constructible ideal with $|G(I)|=r$, $r\geq2$, generated in degree $q$ and let $\Gamma$ be the simplicial complex such that $I=I_{\Gamma}$. We have $I=I_1+I_2$, with $I_1,I_2$  $\mathbf{\underline{1}}-$constructible ideals generated in degree $q$ and $I_1\cap I_2$ is an $\mathbf{\underline{1}}-$constructible ideal generated in degree $q+1$. Let $\Gamma_1,\ \Gamma_2$ be simplicial complexes on the vertex set $[n]$ such that $I_{\Gamma_1}=I_1$ and $I_{\Gamma_2}=I_2$.
	
Since $I_{\Gamma}=I_{\Gamma_1}+I_{\Gamma_2}$ we have that $\Gamma=\Gamma_1\cap\Gamma_2$ and $\Gamma^{\vee}=\Gamma_1^{\vee}\union\Gamma_2 ^{\vee}$. We have to prove that $\Gamma^{\vee}$ is a constructible simplicial complex.

By induction hypothesis, $\Gamma_1^{\vee},\ \Gamma_2^{\vee}$ are constructible simplicial complex of dimension $n-q-1$.

Since $I_1\cap I_2=I_{\Gamma_1}\cap I_{\Gamma_2}=I_{(\Gamma_1^{\vee}\cap\Gamma_2^{\vee})^{\vee}}$, by the induction hypothesis we have that $\Gamma_1^{\vee}\cap\Gamma_2^{\vee}$ is a constructible simplicial complex of dimension $n-q-2$. So $\Gamma^{\vee}=\Gamma_1^{\vee}\union\Gamma_2 ^{\vee}$ is a $(n-q-1)-$dimensional constructible simplicial complex, which ends our proof.
\end{proof}

Now we can complete the diagram for pure simplicial complexes and the connections with the Stanley-Reisner ideal of the Alexander dual:\[\begin{array}{ccccc}
	\Delta\ \mbox{is shellable}&\Longrightarrow&\Delta\ \mbox{is constructible}&\Longrightarrow&\Delta\ \mbox{is Cohen-Macaulay}\\
	\Updownarrow&&\Updownarrow&&\Updownarrow\\
I_{\Delta^{\vee}}\	\mbox{is an ideal with}&\Longrightarrow&I_{\Delta^{\vee}}\ \mbox{is a constructible}&\Longrightarrow&I_{\Delta^{\vee}}\ \mbox{is an ideal with a}\\
	\mbox{linear quotients}&&\mbox{ideal}&&\mbox{linear resolution} 
\end{array}.
\]
\section{Properties of constructible ideals}

In this section we prove that every constructible ideal has a linear resolution and we compute the Betti numbers of a constructible ideal.

For this, we shall need the following two lemmas. 

\begin{Lemma}\cite{Ro}
 Let $R$ be a standard graded $k$-algebra and
	\[0\rightarrow M'\rightarrow M\rightarrow M''\rightarrow0
\]
be an exact sequence of $\mathbb{Z}$-graded $R$-modules. If $M'$ and $M''$ have $q-$linear resolutions, then $M$ has a $q-$linear resolution.
\end{Lemma}

\begin{Lemma}
 Let $R$ be a standard graded $k$-algebra and
	\[0\rightarrow M'\rightarrow M\rightarrow M''\rightarrow0
\]
be an exact sequence of $\mathbb{Z}$-graded $R$-modules. If $M'$ has a $(q+1)-$linear resolution and $M$ has a $q-$linear resolution, then $M''$ has a $q-$linear resolution.
\end{Lemma}
\begin{proof} The exact sequence 
	\[0\rightarrow M'\rightarrow M\rightarrow M''\rightarrow0
\]
yields the exact sequence
	\[\ldots \rightarrow \Tor_{i}^{R}(M',k)_{i+j}\rightarrow  \Tor_{i}^{R}(M,k)_{i+j}\rightarrow \Tor_{i}^{R}(M'',k)_{i+j}\rightarrow \]
	\[\rightarrow \Tor_{i-1}^{R}(M',k)_{i+j}\rightarrow \Tor_{i-1}^{R}(M,k)_{i+j}\rightarrow\ldots
\]
Since $M$ has a $q-$linear resolution, $\Tor_{i}^{R}(M,k)_{i+j}=0$ for all $j\neq q$. For the same reason, $\Tor_{i}^{R}(M',k)_{i+j}=0$ for all $j\neq q+1$.

We get that $\Tor_{i}^{R}(M'',k)_{i+j}=0$ for all $j\neq q$, thus $M''$ has a $q-$linear resolution.
\end{proof}

\begin{Theorem}
\it Let $S=k[x_1,\ldots, x_n]$ be a polynomial ring over a field $k$ and $I$ be a constructible ideal of $S$ generated in degree $q$. Then $I$ has a $q-$linear resolution. 
\end{Theorem}

\begin{proof}Let $\mathbf{\underline{a}}\in\mathbb{Z}_{>0}^n$ such that $I$ is $\mathbf{\underline{a}}-$constructible. We use descending induction on the degree of monomials from the minimal system of generators of the monomial ideal $I$.

If $I=(x_1^{a_1}\ldots x_n^{a_n})$, $I$ has an $|\mathbf{\underline{a}}|-$linear resolution.

Let $q<|\mathbf{\underline{a}}|$ and $I\in\mathcal{J}_{\mathbf{\underline{a}}}$ be an $\mathbf{\underline{a}}-$constructible ideal generated in degree $q$. Now we use induction on the number of monomials from the minimal system of generators.

If $u\in\mathcal{M}_{\mathbf{\underline{a}}}$, $\deg(u)=q$ and $I=(u)$, then $I$ has a $q-$linear resolution. 

Let $I$ be an $\mathbf{\underline{a}}-$constructible ideal generated in degree $q$ with $G(I)=\{u_1,\ldots,u_r\}$, $r\geq2$. There exist $\mathbf{\underline{a}}-$constructible ideals $I_1$ and $I_2$, generated in degree $q$, such that $I=I_1+I_2$ and $I_1\cap I_2$ is an $\mathbf{\underline{a}}-$constructible ideal generated in degree $q+1$. By induction hypothesis, $I_1$ and $I_2$ have $q-$linear resolutions and $I_1\cap I_2$ has a $(q+1)-$linear resolution. From the exact sequence
 \[0\rightarrow I_1\rightarrow I_1\dirsum I_2\rightarrow I_2\rightarrow0
\]
we have, by Lemma $3.1$, that $I_1\dirsum I_2$ has a $q-$linear resolution and from the exact sequence:
\[0\rightarrow I_1\cap I_2\rightarrow I_1\dirsum I_2\rightarrow I_1+I_2\rightarrow0,
\]
$I_1+I_2$ has a $q-$linear resolution, by Lemma $3.2$. So $I$ has a $q-$linear resolution.
\end{proof}

\begin{Corollary}
\it Let $I$ be a constructible ideal generated in degree $q$ and let $I_1$ and $I_2$ be constructible ideals generated in degree $q$ such that $I_1\cap I_2$ is a constructible ideal generated in degree $q+1$ and $I=I_1+I_2$. Then
	\[\beta_i(I)=\beta_i(I_1)+\beta_i(I_2)+\beta_{i-1}(I_1\cap I_2).
\]
\end{Corollary}
\begin{proof}
The exact sequence
\[0\rightarrow I_1\cap I_2\rightarrow I_1\dirsum I_2\rightarrow I_1+I_2\rightarrow0
\]
yields the exact sequence
	\[0\rightarrow \Tor_{i}(I_1\dirsum I_2,k)_{i+q}\rightarrow \Tor_{i}(I_1+ I_2,k)_{i+q}\rightarrow  \Tor_{i-1}(I_1\cap I_2,k)_{i+q}\rightarrow  0.\]
From this sequence we get:
	\[\beta_i(I)=\beta_i(I_1\dirsum I_2)+\beta_{i-1}(I_1\cap I_2).
\]
and, next,
\[\beta_i(I)=\beta_i(I_1)+\beta_i(I_2)+\beta_{i-1}(I_1\cap I_2).
\]
\end{proof}
\section{Polarization of constructible ideals}

We prove that the polarization of a constructible ideal is a square-free constructible ideal.

In the polarization process, homological properties of a monomial ideal are preserved. Since the polarization of a monomial ideal is a square-free monomial ideal, we can apply specific techniques suited for these classes of ideals. 

First we recall the notion of polarization of a monomial ideal and some concepts related to it, following \cite{Ali}.

Let $S=k[x_1,\ldots,x_n]$ be a polynomial ring in $n$ variables over a field $k$ and $u=x_1^{\alpha_1}\ldots x_n^{\alpha_n}$ be a monomial of $S$. The \textit{polarization} of $u$ is the monomial
	\[u^p=\prod_{i=1}^{n}\prod_{j=1}^{\alpha_i}x_{ij}
\]
where $u^p\in k[x_{11},\ldots,x_{1\alpha_1},\ldots, x_{n1},\ldots,x_{n\alpha_n}]$. 

Let $I$ be a monomial ideal of $S$ and $u_1,\ldots, u_m$ be a system of monomial generators for $I$. Then, the ideal generated by monomials $u_1^p,\ldots, u_m^p$ is called \textit{a polarization} of $I$. Since the polarization seems to depend on the system of generators, we consider another system of monomial generators, $v_1,\ldots, v_k$, for the monomial ideal $I$. Let $S'$ be a polynomial ring with suficiently many variables such that for all $1\leq i\leq m$ and all $1\leq j\leq k$, $u_i^p$ and $v_j^p$ are monomials in $S'$. Then $(u_1^p,\ldots,u_m^p)=(v_1^p,\ldots,v_k^p)$ in $S'$. So we note that, in a common polynomial ring extension, all polarizations of a monomial ideal are the same. It follows that we can denote any polarization of a monomial ideal $I$ by $I^p$. If $I$ and $J$ are two monomial ideals of $S$, we write $I^p=J^p$ if a polarization of $I$ and a polarization of $J$ coincide in a common polynomial ring extension.

\begin{Proposition} Let $I$ be a constructible ideal of $S$. Then $I^p$ is a square-free constructible ideal.
\end{Proposition}

\begin{proof} Let $\mathbf{\underline{a}}\in\mathbb{Z}_{>0}^n$ such that $I$ is $\mathbf{\underline{a}}-$constructible.  We use descending induction on the degree of monomials from the minimal system of generators of the monomial ideal $I$.

If $I=(x_1^{a_1}\ldots x_n^{a_n})$, then $I^p$ is a principal square-free monomial ideal, hence is an $\mathbf{\underline{1}}-$constructible ideal.

Let $I\in\mathcal{J}_{\mathbf{\underline{a}}}$ be an $\mathbf{\underline{a}}-$constructible ideal generated in degree $q<|\mathbf{\underline{a}}|$. We use induction on the number of monomials from the minimal system of generators of the ideal $I$. 

If $u\in\mathcal{M}_{\mathbf{\underline{a}}}$, $\deg(u)=q$ and $I=(u)$, the statement is obvious. 

Let $I\in\mathcal{J}_{\mathbf{\underline{a}}}$ be an $\mathbf{\underline{a}}-$constructible ideal generated in degree $q$ with $|G(I)|=r$, $r\geq 2$. There exist $I_1,\ I_2\in\mathcal{J}_{\mathbf{\underline{a}}}$ $\ \mathbf{\underline{a}}-$constructible ideals generated in degree $q$ such that $I=I_1+I_2$ and $I_1\cap I_2$ is an $\mathbf{\underline{a}}-$constructible ideal generated in degree $q+1$. By induction hypothesis $I_1^p,\ I_2^p$ and $I_1^p\cap I_2^p=(I_1\cap I_2)^p$ are $\mathbf{\underline{1}}-$constructible ideals. Hence, $I^p$ is an $\mathbf{\underline{1}}-$constructible ideal.
\end{proof}

\section{Ideals with linear quotients}

In this section we describe the relation between monomial ideals with linear quotients and constructible ideals.

\begin{Proposition}
\it Let $I$ be a monomial ideal of $S$ with linear quotients. Then $I$ is a constructible ideal. 
\end{Proposition}

\begin{proof} We prove by induction on the number of monomials in the minimal system of generators.

If $u$ is a monomial in $S$ and $I=(u)$, then $I$ is a constructible ideal, by definition.

Let $I$ be a monomial ideal, $G(I)=\{u_1,\ldots,u_r\}$, $r\geq2$, be its minimal system of generators and assume that $I$ has linear quotients with respect to the sequence $u_1,\ldots, u_r$. Denote $I_1=(u_1,\ldots, u_{r-1})$ and $I_2=(u_r)$. $I_1,\ I_2$ are constructible ideals, by induction hypothesis.
	\[I_1\cap I_2=\left(u_i u_r/[u_i,u_r]:\ 1\leq i\leq r-1\right)=(x_{l_1}u_r,\ldots,x_{l_t}u_r),
\]
for some $l_i\in\{1,\ldots,n\},\ 1\leq i\leq t,\ t\leq r-1 $. The last equality holds since the ideal $I$ has linear quotients. So $I_1\cap I_2$ is a monomial ideal with linear quotients with at most $r-1$ monomials in the minimal system of generators. By the induction hypothesis, we have that $I_1\cap I_2$ is a constructible ideal and then $I$ is a constructible ideal.
\end{proof}

In \cite{HeTa}, the Betti numbers of an ideal with linear quotients are computed. Namely, for an ideal $I$ with linear quotients, with respect to the sequence $f_1,\ldots, f_s$, we denote by $I_k=(f_1,\ldots, f_k)$ and $L_k=(f_1,\ldots, f_{k-1}):(f_k)$. If $r_k$ is the number of generators of the ideal $L_k$, then
\[\beta_i(I)=\sum_{k=2}^m\left(
\begin{array}{c}
	r_k\\
	i
\end{array}
\right).
\]

We may obtain this formula in the next proposition.

\begin{Proposition}
\it
Let $I$ be a monomial ideal of $S$, $G(I)=\{u_1,\ldots, u_m\}$ be its minimal system of generators and assume that $I$ has linear quotients with respect to the sequence $u_1,\ldots, u_m$. Denote $I_{k}=(u_1,\ldots, u_k)$ and let $r_k$ be the number of generators of the monomial ideal $I_{k-1}:(u_k)$. Then
	\[\beta_i(I)=\sum_{k=2}^m\left(
\begin{array}{c}
	r_k\\
	i
\end{array}
\right),
\]
for all $i\geq1$.
\end{Proposition}

\begin{proof}Since $I_k$ is a monomial ideal with linear quotients with respect to the sequence $u_1,\ldots, u_k$, by Proposition $5.1$, $I_k$ is a constructible ideal.

By Corollary $3.4$,
	\[\beta_i(I_k)=\beta_{i}(I_{k-1})+\beta_{i}((u_k))+\beta_{i-1}(I_{k-1}\cap(u_k))\ \ \ \ \ (*)
\]
for all $k\geq2$, $i\geq1$. 

The multiplication by $u_k$ defines an isomorphism between $I_{k-1}:(u_k)$ and $I_{k-1}\cap(u_k)$. Therefore
	\[|G(I_{k-1}\cap(u_k))|=|G(I_{k-1}:(u_k))|=r_k.
\]
Since $I_k$ has linear quotients, the ideal $I_{k-1}:(u_k)$ is generated by a regular sequence of length $r_k$, and then 
	\[\beta_{i-1}(I_{k-1}:(u_k))=\left(\begin{array}{c}
	r_k\\
	i
\end{array}
\right).
\]
Summing in $(*)$ for $k=2,3,\ldots,m$, we get
	\[\beta_{i}(I)=\sum_{k=2}^{m}\beta_{i-1}(I_{k-1}\cap(u_k))=\sum_{k=2}^m\left(\begin{array}{c}
	r_k\\
	i
\end{array}\right).
\]
\end{proof}

\begin{Proposition} Let $I=(u_1,\ldots,u_r)$ be a monomial ideal of $S$. Then $I$ has linear quotients with respect to the sequence $u_1,\ldots, u_r$ if and only if $I^p$ has linear quotients with respect to the sequence $u_1^p,\ldots, u_r^p$.
\end{Proposition}

\begin{proof} Let $2\leq k\leq r$. By $[$ \cite{Ali}, Lemma $3.3]$ $(u_1,\ldots,u_{k-1}):u_k$ is a prime ideal if and only if $(u_1^p,\ldots,u_{k-1}^p):u_k^p$ is a prime ideal. Since any monomial prime ideal is generated by a sequence of variables, the statement follows.
\end{proof}

\section{Examples}
Now we discuss some examples.

\begin{Example}\rm  The following example of constructible and non-shellable simplicial complex is due to Masahiro Hachimori \cite{Ha}.
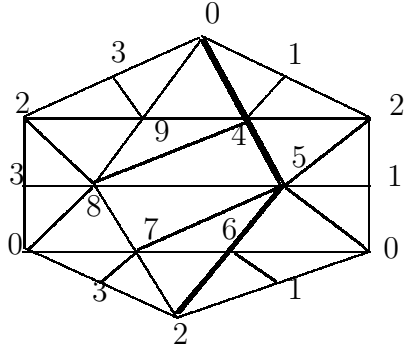
	\[
\]
\unitlength 1mm 
\linethickness{0.4pt}
\ifx\plotpoint\undefined\newsavebox{\plotpoint}\fi 
\begin{picture}(93.75,50.25)(0,0)
\put(45,19.75){\line(1,0){46.25}}
\put(45,28.5){\line(1,0){46.25}}
\multiput(68.5,48.25)(-.074519231,-.033653846){312}{\line(-1,0){.074519231}}
\put(45.25,37.75){\line(0,-1){17.5}}
\multiput(45.25,20.25)(.073636364,-.033636364){275}{\line(1,0){.073636364}}
\multiput(65.5,11)(.098076923,.033653846){260}{\line(1,0){.098076923}}
\put(91,19.75){\line(0,1){18}}
\multiput(91,37.75)(-.069749216,.03369906){319}{\line(-1,0){.069749216}}
\multiput(68.75,48)(-.033687943,-.045508274){423}{\line(0,-1){.045508274}}
\multiput(54.5,28.75)(.033639144,-.054281346){327}{\line(0,-1){.054281346}}
\thicklines
\put(45.5,37.5){\line(1,0){45.5}}
\multiput(69.25,48)(.033653846,-.061698718){312}{\line(0,-1){.061698718}}
\multiput(54.5,29)(.08508403,.03361345){238}{\line(1,0){.08508403}}
\multiput(45.5,37.25)(.035714286,-.033673469){245}{\line(1,0){.035714286}}
\multiput(54.25,28.25)(-.034693878,-.033673469){245}{\line(-1,0){.034693878}}
\multiput(60.25,19.75)(-.03991597,-.03361345){119}{\line(-1,0){.03991597}}
\multiput(60.25,20)(.075961538,.033653846){260}{\line(1,0){.075961538}}
\multiput(79.75,28.5)(.043071161,.033707865){267}{\line(1,0){.043071161}}
\multiput(80,28.25)(.04456522,-.03369565){230}{\line(1,0){.04456522}}
\multiput(90.25,20.5)(.032609,-.032609){23}{\line(0,-1){.032609}}
\multiput(73.25,19.5)(.04910714,-.03348214){112}{\line(1,0){.04910714}}
\multiput(68.75,48)(.03369906,-.061912226){319}{\line(0,-1){.061912226}}
\multiput(79.5,28.25)(-.033687943,-.039598109){423}{\line(0,-1){.039598109}}
\multiput(65.25,11.5)(.033687943,.039598109){423}{\line(0,1){.039598109}}
\multiput(79.5,28.25)(-.03370098,-.041053922){408}{\line(0,-1){.041053922}}
\multiput(69.25,48)(.033670034,-.063131313){297}{\line(0,-1){.063131313}}
\put(69.25,50.25){$0$}
\put(56.75,45){$3$}
\put(44,38.25){$2$}
\put(43.25,28.75){$3$}
\put(43,19.5){$0$}
\put(65,7.5){$2$}
\put(54.25,13.25){$3$}
\put(80.25,13.5){$1$}
\put(93,19){$0$}
\put(93.5,28.5){$1$}
\put(93.75,38){$2$}
\put(80.25,44.75){$1$}
\put(62.5,34.5){$9$}
\put(53.5,24.5){$8$}
\put(72.75,34){$4$}
\put(80.75,31){$5$}
\put(61,21.5){$7$}
\put(71.5,21.5){$6$}
\thinlines
\multiput(75,37.75)(.03355705,.03691275){149}{\line(0,1){.03691275}}
\multiput(60.75,37.75)(-.03348214,.046875){112}{\line(0,1){.046875}}
\thicklines
\multiput(79.5,28.75)(-.033653846,-.041208791){364}{\line(0,-1){.041208791}}
\end{picture}

The simplicial complex is constructible because we can split it by the bold line and we obtain two shellable simplicial complexes $\Delta_1,\ \Delta_2$ of dimension $2$ whose intersection is a shellable $1-$dimensional simplicial complex.
	
The shelling order for the facets of the simplicial complex $\Delta_1$ is:
	\[\{0,3,9\},\{2,3,9\},\{2,8,9\},\{2,3,8\},\{0,3,8\},\{0,7,8\},\{0,3,7\},\{2,3,7\},\{2,6,7\},\]
	\[\{5,6,7\},\{5,7,8\},\{4,5,8\},\{4,8,9\},\{0,4,9\}
\]

For the simplicial complex $\Delta_2$, the shelling order of the facets is
	\[\{0,1,4\},\{1,2,4\},\{2,4,5\},\{1,2,5\},\{0,1,5\},\{0,5,6\},\{0,1,6\},\{1,2,6\}
\]

$\Delta_1\cap\Delta_2$ is the simplicial complex $\langle\{0,4\},\{4,5\},\{5,6\},\{2,6\}\rangle$.

For the Alexander dual of $\Delta_1$, the Stanley-Reisner ideal is
	\[I_{\Delta_1^{\vee}}=(x_1 x_2 x_4 x_5x_6 x_7 x_8,\ x_0x_1x_4x_5x_6x_7x_8,\ x_0x_1x_3x_4x_5x_6x_7,\]\[\ \ \ \ \ \ \ \ \ \ x_0x_1x_4x_5x_6x_7x_9,\ x_1x_2x_4x_5x_6x_7x_9,\ x_1x_2x_3x_4x_5x_6x_9,\]\[\ \ \ \ \ \ \ \ \ \ x_1x_2x_4x_5x_6x_8x_9,\ x_0x_1x_4x_5x_6x_8x_9,\ x_0x_1x_3x_4x_5x_8x_9,\]\[\ \ \ \ \ \ \ \ \ \ x_0x_1x_2x_3x_4x_8x_9,\ x_0x_1x_2x_3x_4x_6x_9,\ x_0x_1x_2x_3x_6x_7x_9,\]\[x_0x_1x_2x_3x_5x_6x_7,\ x_1x_2x_3x_5x_6x_7x_8).\ \ \ \ \ \ \ \ \ \ \ \ 
\]

The Stanley-Reisner ideal for the Alexander dual of $\Delta_2$ is:
\[I_{\Delta_2^{\vee}}=(x_2 x_3 x_5 x_6x_7 x_8 x_9,\ x_0x_3x_5x_6x_7x_8x_9,\ x_0x_1x_3x_6x_7x_8x_9,\]
\[\ \ \ \ \ \ \ \ \ \ x_0x_3x_4x_6x_7x_8x_9,\ x_2x_3x_4x_6x_7x_8x_9,\ x_1x_2x_3x_4x_7x_8x_9,\]
\[x_2x_3x_4x_5x_7x_8x_9,\ x_0x_3x_4x_5x_7x_8x_9).\ \ \ \ \ \ \ \ \ \ \ \ 
\]

The ideals $I_{\Delta_1^{\vee}}$ and $I_{\Delta_2^{\vee}}$ have linear quotients by Theorem $1.6$. 
	\[I_{\Delta_1^{\vee}}\cap I_{\Delta_2^{\vee}}=I_{(\Delta_1\cap\Delta_2)^{\vee}}=(x_1x_2 x_3 x_5 x_6x_7 x_8 x_9,\ x_0x_1x_2 x_3 x_6 x_7x_8 x_9,\]
	\[\ \ \ \ \ \ \ \ \ \ \ \ \ \ \ \ \ \ \ \ \ \ \ \ \ \ \ \ \ \ \ \ \ \ \ \ \ x_0x_1x_2 x_3 x_4x_7 x_8 x_9,\ x_0x_1x_3 x_4 x_5x_7 x_8 x_9).
\]
$I_{\Delta_1^{\vee}}\cap I_{\Delta_2^{\vee}}$ has linear quotients by Theorem $1.6$. The ideal $I_{\Delta^{\vee}}=I_{\Delta_1^{\vee}}+ I_{\Delta_2^{\vee}}$ is a square-free constructible ideal and $I_{\Delta^{\vee}}$ does not have linear quotients since $\Delta$ is not shellable.
	
\end{Example}

\begin{Example}\rm
We consider now the Dunce Hat. It is known that the Dunce Hat is Cohen-Macaulay, but it is not constructible \cite{Ha}.
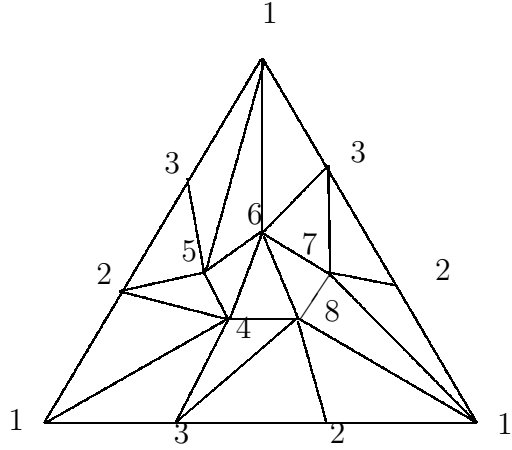
	\[
\]
\unitlength 1mm 
\linethickness{0.4pt}
\ifx\plotpoint\undefined\newsavebox{\plotpoint}\fi 
\begin{picture}(102,63.5)(0,0)
\multiput(70.5,58.5)(-.0337045721,-.0565650645){853}{\line(0,-1){.0565650645}}
\put(41.75,10.25){\line(1,0){57.5}}
\multiput(99.25,10.25)(-.0337278107,.0573964497){845}{\line(0,1){.0573964497}}
\put(70.75,58.75){\line(0,-1){23}}
\multiput(70.75,35.75)(-.0472561,-.03353659){164}{\line(-1,0){.0472561}}
\multiput(63,30.25)(.03350515,-.06443299){97}{\line(0,-1){.06443299}}
\put(66.25,24){\line(1,0){9.5}}
\put(75.75,24){\line(2,3){4}}
\multiput(79.75,30)(-.05487805,.03353659){164}{\line(-1,0){.05487805}}
\multiput(70.75,35.5)(-.03373016,-.08928571){126}{\line(0,-1){.08928571}}
\multiput(66.5,24.25)(-.059036145,-.03373494){415}{\line(-1,0){.059036145}}
\multiput(71,35)(.03358209,-.08022388){134}{\line(0,-1){.08022388}}
\multiput(75.5,24.25)(-.039156627,-.03373494){415}{\line(-1,0){.039156627}}
\multiput(59.25,10.25)(.03365385,.06610577){208}{\line(0,1){.06610577}}
\multiput(66.25,24)(-.12946429,.03348214){112}{\line(-1,0){.12946429}}
\multiput(51.75,27.75)(.15,.0333333){75}{\line(1,0){.15}}
\multiput(63,30.25)(-.0335821,.1865672){67}{\line(0,1){.1865672}}
\multiput(71,58.25)(-.03369565,-.12065217){230}{\line(0,-1){.12065217}}
\multiput(75.5,23.75)(.03348214,-.11830357){112}{\line(0,-1){.11830357}}
\multiput(75.75,24)(.057598039,-.03370098){408}{\line(1,0){.057598039}}
\multiput(99.25,10.25)(-.033703072,.034129693){586}{\line(0,1){.034129693}}
\multiput(79.5,30.25)(.1730769,-.0336538){52}{\line(1,0){.1730769}}
\multiput(79.75,29.75)(-.03125,1.84375){8}{\line(0,1){1.84375}}
\multiput(79.5,44.5)(-.033730159,-.034722222){252}{\line(0,-1){.034722222}}
\put(70.75,63.5){$1$}
\put(57.75,43.5){$3$}
\put(48.75,28.75){$2$}
\put(37,9.25){$1$}
\put(59,7.5){$3$}
\put(79.75,7.5){$2$}
\put(102,8.75){$1$}
\put(93.75,29.25){$2$}
\put(82.5,45){$3$}
\put(68.75,36.75){$6$}
\put(60,31.75){$5$}
\put(67.25,21.5){$4$}
\put(79,23.75){$8$}
\put(76,32.75){$7$}
\end{picture}

The Stanley-Reisner ideal for the Alexander dual of $\Delta$ is:
	\[I_{\Delta^{\vee}}=(x_3x_5x_6x_7x_8,\ x_3x_4x_5x_6x_8,\ x_3x_4x_5x_6x_7,\ x_2x_5x_6x_7x_8,\]	
	\[\ \ \ \ \ \ \ \ \ \ \ x_2x_4x_6x_7x_8,\ x_2x_4x_5x_7x_8,\ x_2x_3x_4x_7x_8,\ x_2x_3x_4x_5x_6,\] 
	\[\ \ \ \ \ \ \ \ \ \ \ x_1x_4x_6x_7x_8,\ x_1x_4x_5x_6x_8,\ x_1x_4x_5x_6x_7,\ x_1x_3x_6x_7x_8,\]
\[\ \ \ \ \ \ \ \ \ \ \ x_1x_2x_5x_6x_7,\ x_1x_2x_4x_5x_8,\ x_1x_2x_3x_7x_8,\ x_1x_2x_3x_5x_7,\]\[x_1x_2x_3x_4x_5)\ \ \ \ \ \ \ \ \ \ \ \ \ \ \ \ \ \ \ \ \ \ \ \ \ \ \ \ \ \ \ \ \ \ \ \ \ \ \ \ 
\]
and it is not a constructible ideal, but it has a linear resolution.
\end{Example}

\begin{Example}\rm Let $I\in k[x_1,\ldots,x_8]$ be the monomial ideal
\[I=(x_1x_2x_5x_6x_7x_8,\ x_2x_3x_5x_6x_7x_8,\ x_2^2x_3x_5x_6x_7,\ x_2^2x_3x_4x_6x_7,\ x_1x_2^2x_3x_6x_7,\ x_2x_3x_4x_5x_7x_8,\]	
	\[\ \ \ \ \ \ \ \ \ \ x_2^2x_3x_4x_7x_8,\ x_1x_2x_3x_4x_7x_8,\ x_1^2x_3x_4x_7x_8,\ x_1^2x_3x_4x_5x_8,\ x_1x_3x_4x_6x_7x_8,\ x_1x_4x_5x_6x_7x_8,\]
  \[\ \ \ \ \ \ \ \ \ \  x_1^2x_4x_5x_6x_8,\ x_1^2x_2x_4x_5x_8,\ x_1x_2^2x_5x_6x_8,\ x_1x_2^2x_3x_6x_8,\ x_1^2x_2^2x_3x_6,\ x_1^2x_2^2x_5x_6,\ x_1^2x_2x_5x_6x_7,\]
  \[x_1^2x_2x_4x_5x_7,\ x_1^2x_2^2x_4x_5)\ \ \ \ \ \ \ \ \ \ \  \ \ \ \ \ \ \ \ \ \ \  \ \ \ \ \ \ \ \ \ \ \  \ \ \ \ \ \ \ \ \ \ \  \ \ \ \ \ \ \ \ \ \ \  \ \ \ \ \ \ \ \ \ \ \ \]
Then $I=I_1+I_2$, where
\[I_1=(x_1x_2x_5x_6x_7x_8,\ x_2x_3x_5x_6x_7x_8,\ x_2^2x_3x_5x_6x_7,\ x_2^2x_3x_4x_6x_7,\ x_1x_2^2x_3x_6x_7,\ x_2x_3x_4x_5x_7x_8,\]
\[ \ \ \ \ \ \ \ \ \ \ \ \ x_2^2x_3x_4x_7x_8,\ x_1x_2x_3x_4x_7x_8,\ x_1^2x_3x_4x_7x_8,\ x_1^2x_3x_4x_5x_8,\ x_1x_3x_4x_6x_7x_8,\ x_1x_4x_5x_6x_7x_8,\]
  \[x_1^2x_4x_5x_6x_8,\ x_1^2x_2x_4x_5x_8)\ \ \ \ \ \ \ \ \  \ \ \ \ \ \ \ \ \ \ \  \ \ \ \ \ \ \ \ \ \ \  \ \ \ \ \ \ \ \ \ \ \ \ \ \ \ \ \ \ \ \  \ \ \ \ \ \ \ \ \ \] 
  and
  \[I_2=(x_1x_2^2x_5x_6x_8,\ x_1x_2^2x_3x_6x_8,\ x_1^2x_2^2x_3x_6,\ x_1^2x_2^2x_5x_6,\ x_1^2x_2x_5x_6x_7,\ x_1^2x_2x_4x_5x_7,\ x_1^2x_2^2x_4x_5)\]
  with
	\[I_1\cap I_2=(x_1x_2^2x_5x_6x_7x_8,\ x_1^2x_2x_5x_6x_7x_8,\ x_1^2x_2x_4x_5x_7x_8,\ x_1^2x_2^2x_4x_5x_8,\ x_1x_2^2x_3x_6x_7x_8,\]
	\[ x_1^2x_2^2x_3x_6x_7)\ \ \ \ \ \ \ \ \  \ \ \ \ \ \ \ \ \ \ \  \ \ \ \ \ \ \ \ \ \ \  \ \ \ \ \ \ \ \ \ \ \ \ \ \ \ \ \ \ \ \  \ \ \ \ \ \ \ \ \ \ \ \ \ \ 
\]  
Since $I_1,\ I_2$ are monomial ideals with linear quotients generated in degree $6$, and $\ \ \ \ I_1\cap I_2$ is a monomial ideal generated in degree $7$ with linear quotients, $I$ is a constructible ideal.

Let $\Delta=\Delta_1\cup\Delta_2$ be the simplicial complex, presented by Ziegler, with $10$ vertices and $21$ facets \cite{Zi}:
	\[\Delta_1: \{1,2,3,4\},\ \{1,2,4,9\},\ \{1,4,8,9\},\ \{1,5,8,9\},\ \{1,4,5,8\},\ \{1,2,6,9\},\ \{1,5,6,9\},\]
	\[ \{1,2,5,6\},\ \{2,5,6,10\},\ \{2,6,7,10\},\ \{1,2,5,10\},\ \{1,2,3,10\},\ \{2,3,7,10\},\ \{2,3,6,7\}
\]
	\[\Delta_2: \{1,3,4,7\},\ \{1,4,5,7\},\ \{4,5,7,8\},\ \{3,4,7,8\}, \ \{2,3,4,8\},\ \{2,3,6,8\},\ \{3,6,7,8\}.
\]
This simplicial complex is constructible, but non-shellable \cite{Zi}.

The polarization of $I$, $I^p$, in the polynomial ring $k[x_{1},x_{1,1},x_2,x_{2,1},x_3,\ldots,x_8]$, with $x_{1,1}=x_9$ and $x_{2,1}=x_{10}$, is the Stanley-Reisner ideal of the Alexander dual associated to the above simplicial complex. Then, the ideal $I$ does not have linear quotients, by Proposition $5.3$.

\end{Example}

It would be interesting to find a large class of constructible ideals which does not have linear quotients.

\ \\ \\
\noindent
Faculty of Mathematics and Informatics,\\
Ovidius University,\\
Bd. Mamaia $124$, $900527$ Constanta\\
ROMANIA\\
E-mail: olteanuandageorgiana@gmail.com

\end{document}